\documentclass[12pt]{amsart}
\usepackage{color}
\usepackage{amsmath}
\usepackage{amssymb}
\usepackage{amsfonts}
\usepackage[latin1]{inputenc}

\DeclareMathSymbol{\subsetneqq}{\mathbin}{AMSb}{36}

 \numberwithin{equation}{section}

\newcommand{\bea}     {\begin{eqnarray}}
\newcommand{\eea}     {\end{eqnarray}}
\newcommand{\beal}     {\begin{align}}
\newcommand{\eeal}     {\end{align}}
\newcommand{\bsp}     {\begin{split}}
\newcommand{\esp}     {\end{split}}
\newcommand{\blll}     {\begin{array}{lll}}
\newcommand{\brcl}     {\begin{array}{rcl}}
\newcommand{\barr}     {\begin{array}}
\newcommand{\earr}     {\end{array}}
\newcommand{\beet}     {\begin{eqnarray*}}
\newcommand{\eec}     {\end{center}}
\newcommand{\bequ}      {\begin{equation}}
\newcommand{\eequ}      {\end{equation}}
\newcommand{\bet}     {\begin{tabular}}
\newcommand{\eet}     {\end{tabular}}
\newcommand{\btab}     {\begin{table}}
\newcommand{\etab}     {\end{table}}
\newcommand{\ben}     {\begin{enumerate}}
\newcommand{\een}     {\end{enumerate}}
\newcommand{\bec}     {\begin{center}}
\newcommand{\beit}     {\begin{itemize}}
\newcommand{\eeit}     {\end{itemize}}
\newcommand{\ov}{\overline }
\newtheorem{theorem}{Theorem}[section]

\newtheorem{lemma}{Lemma}[section]
\newtheorem{proposition}{Proposition}[section]
\newtheorem{definition}{Definition}[section]
\newtheorem{remark}{Remark}[section]

\newcommand{\edoc}{\end{document}}

\def\bspqh{{\dot{B}}^{s, q}_{p}}

\def\S{{\mathcal S}}
\newcommand{\Ref}[1] {(\ref{#1})}
\newcommand{\adelta}{(-\Delta)^{\alpha}}

\newcommand{\rp}{\mathrm{I\kern-.19em R}^{p}}
\newcommand{\rn}{\mathrm{I\kern-.19em R}^{n}}
\newcommand{\rem}{\mathrm{I\kern-.19em R}^{m}}
\newcommand{\rnn}{\mathrm{I\kern-.19em R}^{n\times n}}
\newcommand{\remm}{\mathrm{I\kern-.19em R}^{m\times m}}
\newcommand{\rmn}{\mathrm{I\kern-.19em R}^{m\times n}}
\newcommand{\rnm}{\mathrm{I\kern-.19em R}^{n\times m}}
\newcommand{\rmp}{\mathrm{I\kern-.19em R}^{m\times p}}
\newcommand{\rpn}{\mathrm{I\kern-.19em R}^{p\times n}}
\newcommand{\rpp}{\mathrm{I\kern-.19em R}^{p\times p}}
\newcommand{\rqq}{\mathrm{I\kern-.19em R}^{q\times q}}
\newcommand{\rpq}{\mathrm{I\kern-.19em R}^{p\times q}}
\newcommand{\rppq}{\mathrm{I\kern-.19em R}^{(p+q)\times (p+q)}}
\newcommand{\rnmu}{\mathrm{I\kern-.19em R}^{(n-1)\times (n-1)}}
\newcommand{\rnmd}{\mathrm{I\kern-.19em R}^{(n-2)\times (n-2)}}

\newcommand{\Q}{{\rm \kern.24em \vrule width.02em height1.4ex depth-.05ex \kern-.26em Q}}
\newcommand{\1}{{\rm 1\kern-.25em 1}}
\newcommand{\R}{{\rm I\kern-.25em R}}
\newcommand{\E}{{\rm I\kern-.25em E}}
\newcommand{\N}{{\rm I\kern-.25em N}}
\renewcommand{\H}{{\rm I\kern-.25em H}}
\newcommand{\Z}{{\rm Z \hskip -2.5mm \hskip 1mm Z    }}
\newcommand{\fin}{\rule{3mm}{0mm} \hfill \rule{2mm}{2mm}}
\newcommand{\ph}      {\varphi}

\begin{document}

\title [Global Existence and persistency of the initial regularity ]{Global Existence of Solutions to the $2D$ subcritical dissipative Quasi-Geostrophic equation
 and persistency of the initial regularity }

\author[ R. May. and E. Zahrouni.]{}

\maketitle
\centerline{\scshape May Ramzi\footnote{\email{ramzi.may@fsb.rnu.tn}}}
\medskip
{\footnotesize
 \centerline{ D\'epartement de Math\'ematiques}
\centerline{Facult\'e des Sciences de Bizerte, Tunisie }
}

\medskip

\centerline{\scshape
Zahrouni Ezzeddine\footnote{\email{ezzeddine.zahrouni@fsm.rnu.tn}}}
\medskip
{\footnotesize
 \centerline{ D\'epartement de Math\'ematiques}
\centerline{Facult\'e des Sciences de Monastir, Tunisie }
  }

 \medskip

\medskip

\begin{abstract}
In this paper, we prove that if the initial data $\theta_0$ and
its Riesz transforms ($\mathcal{R}_1(\theta_0)$ and
$\mathcal{R}_2(\theta_0)$) belong to the space
$(\overline{S(\mathbb{R}^2))}^{B_{\infty }^{1-2\alpha ,\infty }}$,
where $\alpha \in ]1/2,1[$, then the $2D$ Quasi-Geostrophic
equation with dissipation $\alpha$ has a unique global in time
solution $\theta$. Moreover, we show that if in addition $\theta_0
\in X$ for some functional space $X$ such as Lebesgue , Sobolev and
Besov's spaces then the solution $\theta$ belongs to the space
$C([0,+\infty [,X).$
\end{abstract}
\section{Introduction and main results}
\smallskip
In this paper, we are concerned with the initial Value-Problem for the
two-dimensional quasi-geostrophic equation with sub-critical dissipation%
\begin{equation*}
\left\{
\begin{array}{c}
\partial _{t}\theta +\left( -\Delta \right) ^{\alpha }\theta + \nabla.(\theta
u)=0\text{ on }\mathbb{R}_{\ast }^{+}\times \mathbb{R}^{2} \\
\theta (0,x)=\theta _{0}(x),~x\in \mathbb{R}^{2}%
\end{array}%
\right. \tag{QG$_{\alpha}$}
\end{equation*}%
where $\alpha \in ]\frac{1}{2},1[$ is a fixed parameter and $\nabla $
denotes the divergence operator with respect to the space variable $x\in
\mathbb{R}^{2}$. The scalar function $\theta $ represents the potential
temperature. The velocity $u=(u_{1},u_{2})$ is divergence free and
determined from $\theta $ through the Riesz transforms%
\begin{equation*}
u=\mathcal{R}^{\bot }(\theta )\equiv \left( -\mathcal{R}_{2}(\theta ),%
\mathcal{R}_{1}(\theta )\right) .
\end{equation*}%
The non local operator $\left( -\Delta \right) ^{\alpha }$ is defined
through the Fourier transform%
\begin{equation*}
\mathcal{F}(\left( -\Delta \right) ^{\alpha }f)(\xi )=\left\vert \xi
\right\vert ^{2\alpha }\mathcal{F}(f)(\xi )
\end{equation*}%
where $\mathcal{F}(f)$ is the Fourier transform of $f$ defined by:%
\begin{equation*}
\mathcal{F}(f)(\xi )=\hat{f}(\xi )=\int_{\mathbb{R}^{2}}f(x)e^{-i\langle
x,\xi \rangle }dx.
\end{equation*}%
To study the existence of the solutions of the equations $(QG_{\alpha })$ we
will follow the Fujita-Kato method. Thus we convert the equations $%
(QG_{\alpha })$ into the fixed point problem%
\begin{equation}
\label{mild}
\theta (t)=e^{-t\left( -\Delta \right) ^{\alpha }}\theta _{0} + {\mathcal{B}}_{\alpha}%
\left[ \theta ,\theta \right] (t).
\end{equation}%
Here $\left( e^{-t\left( -\Delta \right) ^{\alpha }}\right) _{t>0}$ is the
semi-group defined by:%
\begin{equation*}
\mathcal{F}\left( e^{-t\left( -\Delta \right) ^{\alpha }}f\right) (\xi
)=e^{-t\left\vert \xi \right\vert ^{2\alpha }}\mathcal{F}(f)(\xi )
\end{equation*}%
and $\mathcal{B}_{\alpha}$ is the bi-linear operator given by:
\begin{equation}
\label{opB} {\mathcal{B}_{\alpha}}\left[ \theta _{1},\theta
_{2}\right] (t)=-{\mathcal{L}}_{\alpha }(\theta
_{1}\mathcal{R}^{\bot }(\theta _{2}))
\end{equation}%
where, $\; \text{for }v=(v_{1},v_{2}),$
\begin{equation}
\label{linear}
{\mathcal{L}}_{\alpha
}(v)(t)=\int_{0}^{t}e^{-(t-s)\left( -\Delta \right) ^{\alpha }}\nabla .v ds.
\end{equation}%
In the sequel, we mean by a mild solution on $]0,T[$ to the equations $%
(QG_{\alpha })$ with data $\theta _{0}$ a function $\theta $ belonging to
the space $L_{loc}^{2}([0,T[,F_{2})$ and satisfying in ${\mathcal D}^{\prime
}(]0,T[\times \mathbb{R}^{2})$ the equation \Ref{mild} where $F_{2}$ is the
completion of $S(\mathbb{R}^{2})$ with respect to the norm%
\begin{equation*}
\left\Vert f\right\Vert _{F_{2}}\equiv \sup_{x_{0}\in \mathbb{R}^{2}}\left(
\left\Vert 1_{B(x_{0},1)}f\right\Vert _{2}+\left\Vert 1_{B(x_{0},1)}\mathcal{%
R}^{\bot }(f)\right\Vert _{2}\right) .
\end{equation*}%
One of the main property of the equations $(QG_{a})$ is the following
scaling invariance property: if $\theta $ is a solution of $(QG_{a})$ with
data $\theta _{0}$ then, for any $\lambda >0,$ the function $\theta
_{\lambda }(t,x)\equiv \lambda ^{2\alpha -1}\theta (\lambda ^{2\alpha
}t,\lambda x)$ is a solution of $(QG_{a})$ with data $\theta _{0,\lambda
}(x)\equiv \lambda ^{2\alpha -1}\theta _{0}(\lambda x).$ This leads us to
introduce the following notion of \textit{super-critical space}: A Banach
space $X$ will be called \textit{super-critical space} if $S(\mathbb{R}%
^{2})\hookrightarrow X\hookrightarrow S(\mathbb{R}^{2})$ and there exists a
constant $C_{X}\geq 0$ such that
\begin{equation*}
\forall f\in X,~\sup_{0<\lambda \leq 1}\lambda ^{2\alpha -1}\left\Vert
f(\lambda .)\right\Vert _{X}\leq C_{X}\left\Vert f\right\Vert _{X}.
\end{equation*}%
For instance, the Lebesgue space $L^{p}(\mathbb{R}^{2})$ (respectively, the
Sobolev space $H^{s}(\mathbb{R}^{2})$) is \textit{super-critical space} if $%
p\geq p_{c}\equiv \frac{2}{2\alpha -1}$ (respectively, $s\geq s_{c}\equiv
2-2\alpha $). Moreover, one can easily prove that the Besov space $B_{\infty
}^{1-2\alpha ,\infty }(\mathbb{R}^{2})$ is the greatest \textit{%
super-critical space}. The first purpose of this paper, is to prove the
global existence of smooth solutions of the equations $(QG_{a})$ for initial
data in a \textit{sub-critical space} $\mathbf{\tilde{B}}^{\alpha }$ closed
to the space $B_{\infty }^{1-2\alpha ,\infty }(\mathbb{R}^{2}).$ Our space $%
\mathbf{\tilde{B}}^{\alpha }$ is the completion of $S(\mathbb{R}^{2})$ with
respect to the norm%
\begin{equation*}
\left\Vert f\right\Vert _{\mathbf{\tilde{B}}^{\alpha }}\equiv \left\Vert
f\right\Vert _{B_{\infty }^{1-2\alpha ,\infty }}+\left\Vert \mathcal{R}%
^{\bot }(f)\right\Vert _{B_{\infty }^{1-2\alpha ,\infty }}.
\end{equation*}%
Before setting precisely our global existence result, let us
recall some known results in this direction: In \cite{wu2}, J. Wu
proved that for any initial data $\theta _{0}$ in the space
$L^{p}(\mathbb{R}^{2})$ with $p>p_{c}$ the equations $(QG_{a})$
has a unique and global solution $\theta $ belonging to the space
$L^{\infty }([0,+\infty \lbrack ,L^{p}(\mathbb{R}^{2}))$.
Similarly, P. Constantin and J. Wu \cite{constw} showed the global
existence and
uniqueness for arbitrary initial data in the Sobolev space $H^{s}(\mathbb{R}%
^{2})$ where $s>s_{c}.$ Notice that these results don't cover the limit
cases $p=p_{c}$ and $s=s_{c}.$\\
Our global existence result reads as follows.
\begin{theorem}\label{theo1}
Let $\nu =1-\frac{1}{2\alpha }.$ For any initial data $\theta _{0}\in \mathbf{\tilde{B}}%
^{\alpha }$ the equation $(QG_{\alpha })$ has a unique global solution $%
\theta $ belonging to the space $\cap _{T>0}\mathbf{E}_{T}^{\alpha
},$ where $\mathbf{E}_{T}^{\alpha }$ is the completion of
$C_{c}^{\infty }(]0,T]\times \mathbb{R}^{2})$ with respect to the
norm
\begin{equation*}
\left\Vert v\right\Vert _{\mathbf{E}_{T}^{\alpha }}\equiv \sup_{0<t\leq
T}t^{\nu }\left( \left\Vert v(t)\right\Vert _{\infty }+\left\Vert \mathcal{R}%
^{\bot }(v)(t)\right\Vert _{\infty }\right) .
\end{equation*}%
Moreover,
\begin{equation*}
C([0,+\infty \lbrack ,\mathbf{\tilde{B}}^{\alpha }).
\end{equation*}
\end{theorem}
Our second main result is a persistency theorem that states that the
solution $\theta $ given by the previous theorem keeps its initial
regularity. Precisely, our theorem states as follows.
\begin{theorem}\label{theo2}
Let $X$ be one of the following Banach spaces:
\begin{itemize}
\item $X=L^{p}(\mathbb{R}^{2})$ with $1\leq p\leq \infty .$
\item $X=B_{p}^{s,q}(\mathbb{R}^{2})$ with $s>-1$ and $1\leq p,q\leq \infty
. $
\item $X=\dot{B}_{p}^{s,q}(\mathbb{R}^{2})$ with $s>-1$ and $1\leq p,q\leq
\infty .$
\end{itemize}
Assume $\theta _{0}\in \mathbf{\tilde{B}}^{\alpha }\cap X.$ Then
the mild solution $\theta $ of the equation $(QG_{\alpha })$ given
by Theorem \ref{theo1} belongs to the space $L_{loc}^{\infty
}([0,+\infty \lbrack ,X).$ Moreover,
if $\theta _{0}\in \mathbf{\tilde{B}}^{\alpha }\cap \overline{S(\mathbb{R}%
^{2})}^{X}$ then $\theta $ belongs to $C([0,+\infty \lbrack ,\overline{S(%
\mathbb{R}^{2})}^{X}).$
\end{theorem}

As a consequence of the previous theorems, we have the following
theorem that generalizes the existence results of J. Wu
\cite{wu2} and P. Constantin and J. Wu \cite{constw} recalled
above.

\begin{theorem}\label{theo3}
Let $X$ be the Lebesgue space $L^{p}(\mathbb{R}^{2})$ with $p\geq p_{c}=%
\frac{2}{2\alpha -1}$ or the Soblev space $H^{s}(\mathbb{R}^{2})$ with $%
s\geq s_{c}=2-2\alpha .$ Assume $\theta _{0}\in X.$ Then the equation $%
(QG_{\alpha })$ with initial data $\theta _{0}$ has a unique global mild
solution $\theta $ belonging to the space $C([0,+\infty \lbrack ,X).$
\end{theorem}

The remainder of this paper is as follows : in section $2$ we
recall some definitions and we give some useful Lemmas that will
be used in this paper. In section $3$, we prove Theorem
\ref{theo1}. Section {4} is devoted to the proof of Theorem
\ref{theo2} and in section $4$, we will prove Theorem \ref{theo3}.

\vskip 1.0cm
\section{Preliminaries}
\subsection{Notations}
In this subsection, we introduce some notations that will be used
frequently in this paper:
\begin{enumerate}
\item Let $X$ be a Banach space such that
$S(\mathbb{R}^{2})\hookrightarrow X\hookrightarrow S^{\prime
}(\mathbb{R}^{2}).$ We denote by $X_{\mathcal{R}}$ the space
\[
X_{\mathcal{R}}=\{f\in X;~\mathcal{R}^{\perp }(f)\in X^{2}\}
\]%
endowed with the norm%
\[
\left\Vert f\right\Vert _{X_{\mathcal{R}}}=\left\Vert f\right\Vert
_{X}+\left\Vert \mathcal{R}^{\perp }(f)\right\Vert _{X}.
\]%
We recall that $\mathcal{R}^{\perp }(f)=\left( -\mathcal{R}_{2}f,\mathcal{R}%
_{1}f\right) $ where $\mathcal{R}_{1}$ and $\mathcal{R}_{2}$ are
Riesz transforms.

\item Let $T>0,~r\in \lbrack 1,\infty ]$ and $X$ be a Banach
space. $L_{T}^{r}X$ denotes the space $L^{r}([0,T[,X).$ In particular, $%
L_{T}^{r}L^{p}$ will denote the space
$L^{r}([0,T[,L^{p}(\mathbb{R}^{2})).$

\item Let $X$ be a Banach space,~$T>0$ and $\mu \in
\mathbb{R}^{+}.$ we denote by $L_{\mu }^{\infty }([0,T],X)$ the
space of
functions $f:]0,T]\rightarrow X$ such that%
\[
\left\Vert f\right\Vert _{L_{\mu }^{\infty }([0,T],X)}\equiv
\sup_{0<t\leq
T}t^{\mu }\left\Vert f(t)\right\Vert _{X}<\infty \text{ and }%
\lim_{t\rightarrow 0}t^{\mu }\left\Vert f(t)\right\Vert _{X}=0.
\]%
The sub-space $C_{\mu }^{0}([0,T],X)$ of $L_{\mu }^{\infty
}([0,T],X)$ is defined by
\bea
\nonumber
C_{\mu }^{0}([0,T],X)\equiv L_{\mu }^{\infty }([0,T],X)\cap
C(]0,T],X).
\eea
\item Let $A$ and $B$ be two reals functions. The notation
$A\lesssim B$ means that there exists a constant $C,$ independent
of the effective parameters of $A$ and $B,$ such that $A\leq CB.$
\end{enumerate}
\subsection{Besov spaces}\label{section2} 

The standard definition of Besov spaces passes through the Littlewood-Paley dyadic decomposition \cite{bergh}, \cite{frazier}, and \cite{lemarie}. To this end, we take an arbitrary function $\psi \in \S(\R^2) $ whose Fourier transform $\hat{\psi}$ is such that $\; \hbox{supp}( {\hat{\psi}} ) \subset \{\xi,  \frac{1}{2} \leq |\xi| \leq 2 \},\;
\hbox{and for}
\; \xi \neq 0,\; \sum_{j \in \Z} \hat{\psi}(\frac{\xi}{2^j}) = 1,$ and define $\ph \in \S(\R^2) $ by
$ \hat{\ph}(\xi) = 1 - \sum_{ j \geq 0 } \hat{\psi}(\frac{\xi}{2^j}).$
For  $j \in \Z $, we write $\;\ph_j(x) = 2^{2j} \ph( 2^j x )\;$ and $\psi_j(x) = 2^{2j} \psi(2^j x)$
and we denote  the convolution operators $\;S_j\;$ and $\Delta_j,$ respectively, the convolution operators by
$\ph_j$ and $\psi_j$. \\
\noindent
\begin{definition} ~~~\\
\noindent
Let $\;  1 \leq p, q \leq \infty,\;  s \in \R.\;$\\
\noindent
1. A tempered distribution $\; f\;$ belongs to the (inhomogeneous) Besov space $\; B_p^{s,q} \;$  if and only if
\begin{align*}
|| f ||_{B_p^{s,q}} \equiv
||S_0 f||_p + \left(  \sum_{j > 0 } 2^{jsq} || \Delta_j f  ||_p^q \right)^{\frac 1q} < \infty.
\end{align*}
\noindent 2. The homogeneous Besov space $\; \dot{B}_p^{s,q} \;$
is the space of $f\in \S^{\prime}(\mathbb{R}^{2})/_{\R[X]}$ such
that
\begin{align*}
|| f ||_{\dot{B}_p^{s,q}} \equiv \left(  \sum_{j \in \Z } 2^{jsq}
|| \Delta_j f  ||_p^q \right)^{\frac 1q}<\infty,
\end{align*}
Where $ \;\R[X] \;$ is the space of polynomials \cite{peetre}.
\end{definition}
An equivalent definition more adapted to the Quasi-geostrophic equations involves the semigroup $\left( e^{-t\left( -\Delta \right) ^{\alpha }}\right) _{t>0}.$
\begin{proposition}\label{p}
If   $ s < 0$  and $ q = \infty.\; $ Then
\begin{align}
\label{car01}
& f \in \;\dot{B}_p^{s,\infty}\; \Longleftrightarrow \sup_{t > 0 } t^{\frac{-s}{2\alpha}} || e^{-t \adelta} f||_p <
\infty,
\\
\label{car02}
& f \in \;B_p^{s,\infty}\; \Longleftrightarrow \; \forall T > 0,\quad \sup_{ 0 < t < T } t^{\frac{-s}{2\alpha}} || e^{-t \adelta} f||_p  \leq C_T.
\end{align}
\end{proposition}
{Proof : } This proposition can be easily proved by following the
same lines as in the proof of Theorem 5.3 in \cite{lemarie} in the
case of the heaat Kernel. One can see also the proof of
Proposition 2.1 in \cite{miao}.\fin
\subsection{Intermediate results}
\noindent We shall frequently use the following estimates on the
operator  $\;  e^{-t{\adelta}}.$
\begin{proposition}\label{p2}
For $\; t > 0, \;$ we set $\;{\mathcal K}_{t}\;$ the kernel of $\;
e^{-t{\adelta}}.$ Then for all $r\in [1,\infty]$ we have,
\begin{align}
\label{est1}
||  {\mathcal K}_{t} ||_r   = C_{1r}t^{\sigma_r},\\
\label{est1a}
||\nabla {\mathcal K}_{t}||_r  =   C_{2r}  t^{\sigma_r-\frac{1}{2\alpha}},\\
\label{est1b} ||{\mathcal R}_j \nabla {\mathcal K}_{t} ||_r  =
C_{3r} t^{\sigma_r-\frac{1}{2\alpha}},
\end{align}
where $\sigma_r=\frac{1}{\alpha}(\frac{1}{r}-1)$ and $C_{1r},
C_{2r}$ and $ C_{3r}$ are constants independent of $t.$
\end{proposition}
{\bf Proof : }
For the proof of (\ref{est1}-\ref{est1a})  see \cite{miao}. The estimate \Ref{est1b} can be obtained by following the same argument of the proof of Proposition 11.1 in \cite{lemarie}. \fin\\
Following the work of P.G. Lemarié-Rieusset, we introduce the
notion of shift invariant functional space :
\begin{definition}
A Banach space $\;X\;$ is called {\bf  shift invariant functional space } if
\begin{align*}
\bullet  \quad &{\mathcal S}(\mathbb{R}^{2})\hookrightarrow X\hookrightarrow  {\mathcal S}^{'}(\mathbb{R}^{2}), \\
\bullet \quad  &\forall \varphi \in  {\mathcal S}(\mathbb{R}^{2})~\text{and }f\in X, \quad
\left\Vert \varphi \ast f\right\Vert _{X}\leq C_X \left\Vert \varphi
\right\Vert _{1}\left\Vert f\right\Vert _{X}.
\end{align*}
\end{definition}
\begin{remark}
The Lebesgue spaces and Besov spaces are   shift invariant
functional spaces.
\end{remark}
The proof of Theorem \ref{theo1} requires the following lemmas.
\begin{lemma}\label{lemma0}
Let $\;X\;$ be a shift invariant functional space. If  $\; f \in
X\;$ then \bea \label{X0} \sup_{t > 0 } || e^{-t\adelta} f  ||_X
\leq C_X ||  f ||_X. \eea Moreover, if  $\; f \in \ov{\mathcal
S(\mathbb{R}^{2})}^X\;$ then: $e^{-t\adelta}f \in C(] 0, \infty [,
\; \ov{\mathcal S(\mathbb{R}^{2})}^X)$ and  $ e^{-t\adelta}  f
\rightarrow  f $ in $X$ as $t \;\rightarrow 0^+.$
\end{lemma}
{\bf Proof : }
One obtain easily \Ref{X0} from \Ref{est1}.  Let us prove the last
statement. For $ t >0,\;$ we denote by $\mathcal{K}_{t}$ the
kernel of the operator $e^{-t\adelta}.$ Then
$\mathcal{K}_{t}(.)=t^{-\frac{1}{\alpha
}}\mathcal{K}(t^{-\frac{1}{2\alpha }}.)$ where
$\mathcal{K}=\mathcal{K}_{t=1}.$   Since $\mathcal{K}\in
L^{1}(\mathbb{R}^{2})$ and $\int \mathcal{K}(x)dx=1,$ there exists
a sequence $\left( \mathcal{K}_{(n)}\right) _{n}\in \left(
C_{c}^{\infty }(\mathbb{R}^{2})\right) ^{N}$ such that for all
$n,~\int \mathcal{K}_{(n)}(x)dx=1$ and $\left(
\mathcal{K}_{(n)}\right) _{n}\rightarrow \mathcal{K}$ in
$L^{1}(\mathbb{R}^{2}).$ Let $(f_{n})_{n}$
be a sequence in $C_{c}^{\infty }(\mathbb{R}^{2})$ satisfying $%
(f_{n})_{n}\rightarrow f$ in $X.$ Now we consider the functions
$\left(
u_{n}\right) _{n}$ and $u$ defined on $\mathbb{R}^{+\ast }\times \mathbb{R%
}^{2}$ by%
\[
u(t,x)= \mathcal{K}_{t}\ast f \quad \text{ and }\quad u_{n}(t,x)=\mathcal{K}_{(n),t}\ast
f_{n}
\]%
where $\mathcal{K}_{(n),t}(.)=t^{-\frac{1}{\alpha }}\mathcal{K}_{(n)}(t^{-%
\frac{1}{2\alpha }}.)$ and $\ast $ denotes the convolution in $\mathbb{R}%
^{2}.$

One can easily verify that for all $n,$ the function $\hat{u}_{n}(t,\xi )=%
\mathcal{\hat{K}}_{(n)}(t^{\frac{1}{2\alpha }}\xi )\hat{f}_{n}(\xi )$
belongs to the space $C(\mathbb{R}^{+\ast },S(\mathbb{R}^{2}))$ and
satisfies $\hat{u}_{n}(t,.)\rightarrow \hat{f}_{n}$ in $S(\mathbb{R}^{2})$
as $t$ goes to $0^{+}.$ This implies that for all $n,$ $u_{n}$ can be
extended to a function in $C(\mathbb{R}^{+},S(\mathbb{R}^{2}))$ with
$f_{n}$ as value at $t=0.$ Consequently, to conclude the proof of the Lemma,
we just need to show that the sequence $\left( u_{n}\right) _{n}$ converges
to $u$ in the space $L^{\infty }\left( \mathbb{R}^{+},X\right) $. To do
this, we notice that for any $t>0$ and any $n\in \mathbb{N}$ we have,%
\[
u_{n}(t)-u(t)=\mathcal{K}_{(n),t}\ast (f_{n}-f)+\left( \mathcal{K}_{(n),t}-%
\mathcal{K}_{t}\right) \ast f.
\]%
Hence,%
\begin{eqnarray*}
\left\Vert u_{n}(t)-u(t)\right\Vert _{X} &\leq &\left\Vert \mathcal{K}%
_{(n),t}\right\Vert _{1}\left\Vert f_{n}-f\right\Vert _{X}+\left\Vert
\mathcal{K}_{(n),t}-\mathcal{K}_{t}\right\Vert _{1}\left\Vert f\right\Vert
_{X} \\
&\leq& C \left\Vert f_{n}-f\right\Vert _{X}+\left\Vert \mathcal{K}_{(n)}-\mathcal{K%
}\right\Vert _{1}\left\Vert f\right\Vert _{X},
\end{eqnarray*}%
which leads to the desired result. \fin \\
The next lemma will be useful in the sequel.
\begin{lemma}\label{lemma1}
Let $X$ be a shift invariant functional space, $T>0$ and $\mu \geq 0.$ Then, for all $%
f\in L_{\mu }^{\infty }([0,T],X),$ the function
$\mathcal{L}_{\alpha }(f)$
belongs to $L_{\mu ^{\prime }}^{\infty }([0,T],X_{\mathcal{R}})$ and satisfies%
\[ \left\Vert \mathcal{L}_{\alpha }(f)\right\Vert _{L_{\mu
^{\prime }}^{\infty }([0,T],X_{\mathcal{R}})}\leq C\left\Vert
f\right\Vert _{L_{\mu }^{\infty }([0,T],X)}
\]%
where $\mu ^{\prime }=\mu -1+\frac{1}{2\alpha }$ and $C$ is a
constant depending only on $\alpha $ and $X.$ Moreover, if $f$
belongs to $L_{\mu
}^{\infty }([0,T],\overline{S(\mathbb{R}^{2})}^{X})$ then $\mathcal{L}%
_{\alpha }(f)$ belongs to the space $C_{\mu ^{\prime }}^{0}([0,T],(\overline{%
S(\mathbb{R}^{2})}^{X})_{\mathcal{R}}).$
\end{lemma}
{\bf Proof : } The first assertion is a an immediate consequence
of estimates (\ref{est1a})-(\ref{est1b}). The last assertion can
be easily proved by using the previous lemma and the Lebesgue's
dominated convergence theorem, we
left details to the reader. \fin 
\begin{lemma}\label{lemma2} Let $\;T > 0.\;$ \\ The following
assertions hold true:
\begin{enumerate}
\item The linear operator $e^{-t\adelta}$ is continuous from $\;\mathbf{\tilde{B}}^{\alpha }\;$ to $\;\E_T^{\nu} .\;$ 
\item The bilinear operator $\mathcal{B}_{\alpha}$ is continuous
from $\;\E_T^{\nu} \times \E_T^{\nu} \rightarrow \;\E_T^{\nu}\;$
and its norm is independent of $T.$
\end{enumerate}
\end{lemma}
{\bf Proof :} The first assertion follows from the
characterization of Besov spaces by the  kernel
$\;e^{-t\adelta}\;$ and the definition of
$\;\mathbf{\tilde{B}}^{\alpha }\;$
The second assertion, is a direct consequence of the previous lemma and the fact that $\;\E_T^{\nu}=C_{\nu }^{0}([0,T],\left( C_{0}(\mathbb{R}^{2})\right) _{\mathcal{R}})$  \fin\\
The following Lemma, which is a direct consequence of the
preceding one will be useful in the proof of Theorem \ref{theo2}.
\begin{lemma}\label{lemma3} Let $\; \theta_{0} \in \tilde{B}^{\alpha}. \;$ The sequence $\phi_n(\theta_{0})$ defined by
\begin{eqnarray*}
\begin{array}{rcl}
\phi_0(\theta_{0}) & = &\;e^{-t\adelta}\theta_{0},\;\\
\phi_{n+1}(\theta_{0}) & = &\;e^{-t\adelta}\theta_{0} + {\mathcal
B}_{\alpha}[\phi_n(\theta_{0}),\phi_n(\theta_{0})],
\end{array}
\end{eqnarray*}
belongs to $\; \bigcap_{T > 0} \E_T^{\nu}. \;$ Moreover, there exists a constant $\; \mu_0 > 0 \;$  (depending only on  $\alpha$ ) such that if for some $T > 0$ we have $\; ||\phi_0(f)||_{\E_T^{\nu}} \leq \mu_0\;$ then $\; \forall n \in \N^{*},$
\begin{eqnarray}
\label{suite1}
||\phi_n(\theta_{0})||_{\E_T^{\nu}}  & \leq &  2 ||\phi_0(\theta_{0})||_{\E_T^{\nu}},\\
\label{suite2} ||\phi_{n+1}(\theta_{0}) -
\phi_{n}(\theta_{0})||_{\E_T^{\nu}} & \leq & \frac{1}{2^n}.
\end{eqnarray}
In particular, the sequence $(\phi_n(\theta_{0}))_n $ converges in
the space $\E_T^{\nu}$ and its limit $\theta$ is a mild solution
to the equation $(QG_{\alpha })$ with initial data $\theta_{0}$.
\end{lemma}
\noindent The following elementary lemma will play a crucial role
in this paper.
\begin{lemma}\label{lemma5} ( {\bf Gronwall type Lemma })~~ Let $T > 0,\; c_1, c_2 \geq 0,\; \kappa \in ] 0, 1 [$ and
$\; f\in  L^{\infty}(0,T) \; $  such that for all $\; t \in [ 0 , T ] \; $
\begin{equation}
f(t) \leq c_1 + c_2 \int_0^t \frac{f(s)}{(t-s)^{\kappa}}.
\end{equation}
Then
\begin{equation}
\label{ineg_lem5}
\;\forall t \in [ 0 , T ], \quad f(t) \leq 2 c_1 e^{\nu t},
\end{equation}
where $\; \nu = \nu_{\kappa,c_2} > 0.\;$
\end{lemma}
{\bf Proof :}
Let $\nu > 0$ to be precise in the sequel and consider the function $\;g\;$ defined on $\; [ 0 , T ]\;$  by
$$
g(t) = \sup_{0 < s < t } e^{-\nu s } f(s).
$$
Clearly, we have
$$
\begin{array}{rcl}
g(t) &\leq &  c_1 + c_2 \int_0^t \frac{e^{-\nu(t-s)}}{(t-s)^{\kappa}} g(s) ds,\\
     & \leq & c_1 + c_2 \gamma_{\kappa} \nu^{\kappa-1} g(t),
\end{array}
$$
where  $\;\gamma_{\kappa} \; = \; \int_0^{\infty}
\frac{e^{-t}}{t^{\kappa}}.$ Thus, if we choose $\nu >0$ such that
$ c_2 \gamma_{\kappa} \nu^{\kappa-1}  = \frac 12\;$, we get the estimate (\ref{ineg_lem5}). \fin \\
\vskip 2.0mm
\begin{lemma}\label{lemme6} {(\bf Maximum Principal)} ~~~~\\
Let $\theta $ be a mild solution to the equation (\ref{mild})
belonging to the space $C([0,T],\left(
C_{0}(\mathbb{R}^{2})\right) _{\mathcal{R}})$. Then $\; \forall  t
\in [ 0,  T ], \;$ we have  \bequ \label{max1}
||\theta(t)||_{\infty} \leq ||\theta_0||_{\infty}, \eequ \bequ
\label{max2}
 ||{\mathcal{R}^{\bot }}(\theta)(t)||_{\infty} \leq  2||{\mathcal{R}^{\bot }}( \theta_0)||_{\infty} e^{\eta t},
\eequ where $\; \eta  = \eta_{\alpha, ||\theta_0||_{\infty}} > 0.
\;$
\end{lemma}
{\bf Proof : } The inequality \Ref{max1} is proved in
\cite{resnick}, \cite{const2001} and \cite{wu2}, for sufficiently
smooth solution $\theta .$  To prove it in our case, we will
proceed by linearization of the equations and regularization of
the initial data. We consider a sequence of {\bf linear
system } $(QGL_{n})_{n}:$%
\begin{equation}
\left\{
\begin{array}{c}
\partial _{t}v-\adelta v+\nabla .(u_{n}v)=0 \\
v(0,.)=\theta _{n}(.).%
\end{array}%
\right.  \tag{QGL$_{n}$}
\end{equation}%
where $\left( \theta _{n}\right) _{n}$ is a given sequence in $C_{c}^{\infty }(%
\mathbb{R}^{2})$ converging to $\theta (0)$ in the space
$L^{\infty }\left( \mathbb{R}^{2}\right) $ and $u_{n}=\omega
_{n}\ast \mathcal{R}^{\bot }(\theta )$
with $\omega _{n}(.)=n^{2}\omega (n.)$ where $\omega \in C_{c}^{\infty }(%
\mathbb{R}^{2})$ and $\int \omega dx=1.$

Let $n\in \mathbb{N}.$ By converting the system  (QGL$_{n}$ ) into the integral
equation
\begin{equation}
v(t)=e^{-t\adelta}\theta _{n}-\int_{0}^{t}\nabla
.e^{-(t-s)\adelta}(u_{n}v)ds \tag{IQGL$_{n}$}
\end{equation}%
and by following a standard method, one can easily prove that the system   (QGL$_{n}$ )
 has a unique global solution $v_{n}\in \cap _{k\in
\mathbb{N}}C^{\infty }([0,T],H^{k}(\mathbb{R}^{2})).$ Hence we are allowed
to make the following computations: Let $p\in \lbrack 2,\infty \lbrack .$
For any $t\in \lbrack 0,T]$ we have%
\begin{eqnarray*}
\frac{1}{p}\frac{d}{dt}\left\Vert v_{n}(t)\right\Vert ^{p}
&=&-\int \left( \adelta v\right) v\left\vert v\right\vert
^{p-2}dx-\int \nabla
.(u_{n}v)v\left\vert v\right\vert ^{p-2}dx \\
&\equiv &I_{1}(t)+I_{2}(t).
\end{eqnarray*}%
Firstly, a simple integration by parts implies that $I_{2}(t)=-I_{2}(t)$ and
so
\[
I_{2}(t)=0.
\]%
Secondly, by the positivity Lemma  (see \cite{resnick} and
\cite{cordoba2004}), we have
\[
I_{1}(t)\leq 0.
\]%
Therefore,%
\[
\sup_{t\in \lbrack 0,T]}\left\Vert v_{n}(t)\right\Vert _{p}\leq \left\Vert
\theta _{n}\right\Vert _{p}.
\]%
Letting $p\rightarrow +\infty ,$ yields%
\[
\sup_{t\in \lbrack 0,T]}\left\Vert v_{n}(t)\right\Vert _{\infty }\leq
\left\Vert \theta _{n}\right\Vert _{\infty }.
\]%
Consequently, to obtain the inequality \Ref{max1} we just need to
show that the sequence $\left( v_{n}\right) _{n}$ converges to the
function $\theta $ in the space $L^{\infty }([0,T],L^{\infty
}(\mathbb{R}^{2})).$ To do this, we consider the sequence $\left(
w_{n}\right) _{n}=\left( v_{n}-\theta \right)
_{n}.$ Let $t\in \lbrack 0,T]$ and $n\in \mathbb{N}$. We have%
\begin{eqnarray*}
w_{n}(t) &=&e^{-t\adelta}\left( w_{n}(0)\right)
-\int_{0}^{t}\nabla .e^{-(t-s)\adelta}(\left( u_{n}-\mathcal{R}^{\bot }%
(\theta )\right) v_{n})ds \\
&&-\int_{0}^{t}\nabla .e^{-(t-s)\adelta}(\mathcal{R}^{\bot
}(\theta )w_{n})ds.
\end{eqnarray*}%
Thus, by using the Young inequality and Proposition \ref{p2}, we easily get%
\[
\left\Vert w_{n}(t)\right\Vert _{\infty }\leq \left\Vert \theta _{n}-\theta
(0)\right\Vert _{\infty }+C_{\alpha }T^{\nu }A_{n}B_{n}+C_{\alpha }M_{\theta
}\int_{0}^{t}\frac{\left\Vert w_{n}(s)\right\Vert _{\infty }}{\left(
t-s\right) ^{1/2\alpha }}ds
\]%
where $C_{\alpha }$ is a constant depending only on $\alpha ,$
\[
A_{n}=\sup_{0\leq t\leq T}\left\Vert u_{n}(t)-\mathcal{R}^{\bot
}(\theta )(t)\right\Vert _{\infty },
\]%
\[
B_{n}=\sup_{0\leq t\leq T}\left\Vert v_{n}(t)\right\Vert _{\infty }
\]%
and
\[
M_{\theta }=\sup_{0\leq t\leq T}\left\Vert \mathcal{R}^{\bot
}(\theta )(t)\right\Vert _{\infty }.
\]%
Applying Lemma \ref{lemma5}, we get%
\[
\sup_{0\leq t\leq T}\left\Vert w_{n}(t)\right\Vert _{\infty }\leq C\left[
\left\Vert \theta _{n}-\theta (0)\right\Vert _{\infty }+C_{\alpha }T^{\nu
}A_{n}B_{n}\right]
\]%
where $C$ is a constant depending on $\alpha ,T$ and $\theta $
only.

Therefore, to obtain the desired conclusion, we just have to notice that the
sequence $(B_{n})_{n}$ is bounded and that $A_{n}\rightarrow 0$ as $%
n\rightarrow \infty $ thanks to the uniform continuity of the function $%
\mathcal{R}^{\bot }(\theta )$\ on $[0,T]\times \mathbb{R}^{2},$
which is a consequence of the fact $\mathcal{R}^{\bot }(\theta
)\in C([0,T], C_{0}(\mathbb{R}^{2}))$

Now, let us establish the inequality \Ref{max2}. For any $t\in \lbrack 0,T],$ we have

\[
\mathcal{R}^{\bot }(\theta )(t)=e^{-t\adelta}\left(
\mathcal{R}^{\bot }(\theta
)(0)\right) -\int_{0}^{t}\mathcal{R}^{\bot }\nabla .e^{-(t-s)\adelta}(%
\mathcal{R}^{\bot }(\theta )\theta )ds.
\]%
Applying the Young inequality and \Ref{est1b}, we get%
\bea
\nonumber
\left\Vert \mathcal{R}^{\bot }(\theta )(t)\right\Vert _{\infty
}\leq \left\Vert \mathcal{R}^{\bot }(\theta )(0)\right\Vert
_{\infty }+C\left\Vert
\theta (0)\right\Vert _{\infty }\int_{0}^{t}\frac{\left\Vert \mathcal{R}^{\bot }%
(\theta )(s)\right\Vert _{\infty }}{\left( t-s\right) ^{1/2\alpha
}}ds
\eea
where the constant $C$  depends only on $\alpha$i.
Hence, Lemma \ref{lemma5} leads the desired inequality. \fin \\

\section{Proof of Theorem \ref{theo1}}
~~\\
According to Lemma \ref{lemma2}, there exists $T>0$ such that
$\left\Vert e^{-t\left( -\Delta \right) ^{\alpha }}\theta
_{0}\right\Vert _{\mathbf{E}_{T}^{\nu }}\leq \mu _{0}$ where $\mu
_{0}$ is the real defined by Lemma \ref{lemma3}. Therefore, the
same lemma ensures  that the equation $(QG_{\alpha })$  with
initial data $\theta _{0}$ has a mild solution $\theta $ belonging
to the space $\mathbf{E}_{T}^{\nu }.$ Following a standard
arguments (see for example the proof of the
\cite[Lemma]{lemarie}), the uniqueness of the solution $\theta $
can be easily deduced from the continuity of the operator
$\mathcal{B}_{\alpha }$ on the space $\mathbf{E}_{T}^{\nu }.$
Hence, there exists a unique maximal solution,
$$
\theta \in \bigcap _{0<T<T^{\ast }}\mathbf{E}_{T}^{\nu }.
$$
where $\;T^{*}\;$ is the maximal time existence. Let us show that,
$$
\theta \in C([0, T^{*}), \mathbf{\tilde{B}}^\alpha ).
$$
Thanks to the embedding,
$$
(C_0(\R^2))_{\mathcal{R}} \subset \mathbf{\tilde{B}}^\alpha ,
$$
and Lemma \ref{lemma0}, we just need to prove the continuity of,
$$
N(\theta)(t) = {\mathcal B}_{\alpha }[\theta,\theta](t),
$$
at $t = 0^+$ in the space $\mathbf{\tilde{B}}^\alpha$.  Even more,
we show that
$$
\lim_{t\rightarrow 0^+} N(\theta )(t) = 0,\quad \hbox{in}\quad
\mathbf{\tilde{B}^\alpha }.
$$
For that, we use Proposition \ref{p2}, the Young inequality and
 estimates $ \Ref{est1a}-\Ref{est1b}$, to get
\begin{align*}
|| N(\theta )(t)||_{\mathbf{\tilde{B}}^\alpha }
 \lesssim
  \sup_{ 0 < t^{'} < 1}
 { t^{'}}^{\nu}
  \int_0^t ( t + t^{'} -\tau )^{-\frac{1}{2\alpha }} \tau^{-2\nu} d\tau \quad ||\theta ||_{\E_t^{\nu}}^2,
\end{align*}
and hence we obtain,
\begin{align}
\label{est5aa} || N(\theta)(t)||_{\mathbf{\tilde{B}}^\alpha }
\lesssim ||\theta ||_{\E_t^{\nu}}^2.
\end{align}
Since the right hand side of \Ref{est5aa} goes to $0$ as $t$ goes $0^+$ we obtain the desired result.\\
It remains to show that the solution $\theta $ is global, that is
$\; T^{\ast}=\infty .\;$ We argue by contradiction. If  $T^{\ast
}< \infty $ then, from  Lemma \ref{lemma3}, we must have,
$$
\forall \;  0 < t_{0} < T^{*},\quad \left\Vert
e^{t\adelta}\theta (t_{0})\right\Vert _{\mathbf{E}_{T^{\ast }-t_{0}}^{\nu }}\geq \mu_{0},
$$
 which yields by the Young inequality
\bequ \label{est6b} \left\Vert \theta (t_{0})\right\Vert _{\infty
}+\left\Vert {\mathcal R}^{\bot }(\theta)(t_{0})\right\Vert
_{\infty }\geq \frac{c}{(T^{\ast }-t_{0})^{\nu }}, \eequ where
$c>0$ is a universal constant. Which contradicts the maximum
principal ( Lemma \ref{lemma5}).
\section{Proof of Theorem \ref{theo2}}\label{sectreg}

 Along this section, we consider  $\theta_0$ a given initial data
 belonging to the space $\mathbf{\tilde{B}}^\alpha$ and we denote by $\theta$ the solution to
 the equation (QG$_{\alpha}$) given by Theorem
 \ref{theo1}. We will establish the persistency of the regularity of the initial data. That is, if moreover $\theta_0 \in X $ for a suitable Banach spaces $X$ then  the solution $\theta \in C([0,\infty), X).\;$
\subsection{Propagation of the $\;L^p\;$ regularity}
In this subsection we will prove the propagation of the initial $
  L^p  $ regularity. Precisely, we prove the following proposition.
\begin{proposition}\label{propLp}
Let   $ X  =  L^p \ ; $     with   $ \;p \in [ \; 1, \infty\;].$
If $\; \theta_0 \in  X \;$ then $\; \theta  \in \bigcap_{T >
0}L^{\infty}([ 0, T], X). \; $ Moreover, if  $\; \theta_0 \in
\overline{S(\mathbb{R}^{2})}^X \;$ then   $\; \theta  \in  C([ 0,
\infty ), \overline{S(\mathbb{R}^{2})}^X ) \; $
\end{proposition}
{\bf Proof : }
assume  $\; \theta_0 \in  X \;$ and let $\;  T > 0.\; $ We consider the Banach spaces
 $\;  \mathbf{Z}_1 = \E_{T}^{\nu} \; $ and  $\;   \mathbf{Z}_2 = L^{\infty}([ 0, T], X) \; $  endowed respectively with the norm
$$
|| v ||_{\mathbf{Z}_{1}} = \sup_{0 < t < T} e^{-\lambda t} t^{\nu} || v(t)||_{\infty} \quad \hbox{and} \quad
|| v ||_{\mathbf{Z}_{2}} = \sup_{0 < t < T} e^{-\lambda t}|| v(t)||_{p},
$$
where $\lambda > 0$ to be fixed later. We consider the linear
integral equation, \bequ \label{linpsi} v = \Psi _{\theta }(v)
\equiv e^{t\Lambda ^{2\alpha }}\theta _{0}+ \mathcal{B}_{\alpha
}[\theta ,v]. \eequ  Let $k\in \{1;2\}$. According to Lemma
\ref{lemma1}, the affine functional  $\Psi_{\theta} : \;
\mathbf{Z}_{k} \rightarrow \mathbf{Z}_{k}\;$ is continuous. Let us
estimate the norm of its linear part,
\[
K_{\theta }(v)=\mathcal{B}_{\alpha }[\theta ,v].
\]%
Let $\varepsilon >0$ to be chosen later. A direct computation using \Ref{est1a} gives,%
\begin{eqnarray*}
\left\Vert K_{\theta }\right\Vert _{\mathcal{L}(\mathbf{Z}_{1})}
&=&\sup_{\left\Vert v\right\Vert _{\mathbf{Z}_{1}}}\left\Vert
K_{\theta
}(v)\right\Vert _{\mathbf{Z}_{1}} \\
&\leq &C_{1}\sup_{0<t<T}t^{\nu }\int_{0}^{t}(t-\tau )^{-\frac{1}{2\alpha }%
}\tau ^{-2\nu }e^{-\lambda (t-\tau )}\left\Vert \theta \right\Vert _{\mathbf{%
E}_{\tau }^{\nu }}d\tau  \\
&\leq &C_{2}\left( \left\Vert \theta \right\Vert
_{\mathbf{E}_{\varepsilon
}^{\nu }}\sup_{0<t<\varepsilon }t^{\nu }\int_{0}^{t}(t-\tau )^{-\frac{1}{%
2\alpha }}\tau ^{-2\nu }d\tau +T^{\nu }\varepsilon ^{-2\nu
}\left\Vert
\theta \right\Vert _{\mathbf{E}_{T}^{\nu }}\lambda ^{-\nu }\Gamma (\nu )\right)  \\
&\leq &C_{3}\left( \left\Vert \theta \right\Vert
_{\mathbf{E}_{\varepsilon
}^{\nu }}+T^{\nu }\varepsilon ^{-2\nu }\lambda ^{-\nu }\left\Vert \theta \right\Vert _{%
\mathbf{E}_{T}^{\nu }}\right) .
\end{eqnarray*}%
where the constants $C_{1},C_{2}$ and $C_{3}$ depend only on
$\alpha .$
Similarly, we prove the estimate,
 \[
\left\Vert K_{\theta }\right\Vert
_{\mathcal{L}(\mathbf{Z}_{2})}\leq C\left( \left\Vert \theta
\right\Vert _{\mathbf{E}_{\varepsilon }^{\nu }}+T^{\nu
}\varepsilon ^{-2\nu }\lambda ^{-\nu }\left\Vert \theta \right\Vert _{%
\mathbf{E}_{T}^{\nu }}\right)
\]
where $C$ is a constant depending only on $\alpha .$ Since,
$\;||\theta||_{\E_{\epsilon}^{\nu}} \rightarrow 0\;$ as
$\;\epsilon \rightarrow 0^+,\;$  one can choose, successively,
$\;\epsilon \;$ small enough
 and  $\lambda \;$ large enough so that $\Psi_{\theta}$  becomes a contraction on $\;\mathbf{Z}_{1}\;$ and $\;\; \mathbf{Z}_{2}\;$ and therefore on $\;\mathbf{Z}_{1}\cap \mathbf{Z}_{2}.\;$ Let  $v_1$ and $v_{1,2}$ be the unique fixed point of $\Psi_{\theta}$ respectively in $\;\mathbf{Z}_{1}\;$ and  $\mathbf{Z}_{1}\cap \mathbf{Z}_{2}.\;$ Now, since $\;\mathbf{Z}_{1} \cap \mathbf{Z}_{2} \subset \mathbf{Z}_{1}\;$ then  $v_1 = v_{1,2}.$
 Moreover, by construction $\theta$ is a fixed point of $\Psi_{\theta}$ in $\;\mathbf{Z}_{1}\;$ thus  $\theta = v_1 = v_{1,2}$ and hence $\; \theta \in  L^{\infty}([ 0, T], X). \; $  \\
The proof of the last statement of the proposition is identically
similar, we have only to replace $\;\mathbf{Z}_{2}\;$  by   $\; C([ 0, T] ,\overline{S(\mathbb{R}^2})^X ). \; $ \fin
\subsection{Propagation of $\;\bspqh\;$ regularity for $\; s > 0$}

In this section, we prove an abstract result, which implies in particular the persistence of the $\;\bspqh\;$ regularity for $s>0.$ Our result states as follows :
\begin{proposition}\label{propositionabstraite}
Let $X$ be a shift invariant functional space such that for a constant $C$%
\bequ ~\forall f,g\in X\cap L^{\infty}(\mathbb{R}^{2}),\quad
\left\Vert fg\right\Vert _{X}\leq C\left( \left\Vert f\right\Vert
_{\infty }\left\Vert g\right\Vert _{X}+\left\Vert g\right\Vert
_{\infty }\left\Vert f\right\Vert _{X}\right) . \eequ If the
initial data $\theta _{0}$ belongs to $X_{\mathcal{R}}$ then the
solution $\theta $ belongs to $\bigcap_{T > 0} L^{\infty} ([0,T],
X_{\mathcal R}) $. Moreover, if $\theta _{0}$ belongs to
$(\overline{S(\mathbb{R}^2)}^X)_{\mathcal{R}}$ then $\theta $
belongs to $ C\left( \mathbb{R}^{+},
(\overline{S(\mathbb{R}^2)}^X)_{\mathcal{R}}\right) $).
\end{proposition}
\vskip 5.0 mm
The proof of this proposition relies essentially on the two followings
lemmas. The first one is an elementary compactness lemma :
\begin{lemma}\label{compact}
Let $\lambda >0$ and $K$ a compact subset of
$\mathbf{\tilde{B}}^{\alpha }.$
Then there exists $\delta =\delta (K,\lambda )>0$ such that%
\[
~\forall f\in K,~~
\left\Vert e^{-t\adelta }f\right\Vert _{\E_{\delta }^{\nu }}\leq
\lambda.
\]
\end{lemma}
{\bf Proof :}
For $n\in \mathbb{N}^{\ast },$ we set
\[
V_{n}=\left\{ f\in \mathbf{\tilde{B}}^{\alpha }, ~~\left\Vert
e^{-t\adelta } f\right\Vert _{\E_{1/n}^{\nu }}<\lambda \right\} .
\]%
We claim that, $\; \forall n \in \N^*, \; V_n \; $ is an open
subset of $\mathbf{\tilde{B}}^{\alpha }$ and $\cup _{n}V_{n} =
\mathbf{\tilde{B}}^{\alpha }.$ This follows easily from the
continuity of the linear operator $e^{-t\adelta} $
 from $\mathbf{\tilde{B}}^{\alpha }$ into $\E_{T}^{\nu }$ for all $T>0$ and the
propriety
\[
\forall f \in \mathbf{\tilde{B}}^{\alpha },\quad
\lim_{T\rightarrow 0}\left\Vert e^{-t\adelta} f \right\Vert
_{\E_{T}^{\nu }}=0.
\]%
Thus, since $K$ is a compact subset of $\mathbf{\tilde{B}}^{\alpha
},$ there exists a finite subset $I \subset \mathbb{N}^{\ast }$
such that $K \subset \cup _{I}V_{n} = V_{n^{\ast }}$ where
$n^{\ast }= \max( n \in I).$ Hence, we conclude that the choice
$\delta =1/n^{\ast }$ is suitable. \fin
\vskip 5.0mm The second lemma establishes a local in time
propagation of the $X$ regularity.
\begin{lemma}\label{localpersistence} Let $X$ be as in the Proposition \ref{propositionabstraite}.
If $\theta _{0}$ belongs to $X_{\mathcal{R}}$ (resp.\;
$(\overline{S(\mathbb{R}^2)}^X)_{\mathcal{R}} $) then there exists
$\delta =\delta(X,\alpha )>0$ such that the solution
$\theta \in L^{\infty }\left( \left[ 0,\delta \right] ,X_{\mathcal{R}}\right) $ ( resp. $C\left( %
\left[ 0,\delta \right] ,
(\overline{S(\mathbb{R}^2)}^X)_{\mathcal{R}}\right)$. Moreover,
the time $\; \delta $ is bounded below by,%
\[
\sup \left\{ T>0, \; \left\Vert e^{-t\adelta }\theta _{0}\right\Vert
_{\E_{T}^{\nu }}\leq \mu \right\},
\]%
where $\mu $ is a non negative constant depending on $X$ and $\alpha $
only.
\end{lemma}
{\bf Proof :} let us consider the case of $\theta _{0}\in X_{\mathcal{R}}$.
The proof in the other case is similar. \\
Let $\mu \in \left] 0,\mu _{0}\right[ $ to be chosen later and let $T>0$
such that $\left\Vert e^{-t\adelta}\theta _{0}\right\Vert _{\mathbf{%
E}_{T}^{\nu }}\leq \mu. $ According to the Lemma \ref{lemma3}, the sequence $\left( \phi
_{n}(\theta _{0})\right) _{n}$ converges in $\E_{T}^{\nu }$ to the solution $%
\theta $ and satisfies the following estimates%
\begin{eqnarray}
\label{est20a}
\sup_{n}\left\Vert \phi _{n}(\theta _{0})\right\Vert _{\mathbf{E}_{T}^{\nu
}} &\leq &\mu \\
\label{est20b}
\forall n \in \N,\quad \left\Vert \phi _{n+1}(\theta _{0})-\phi _{n}(\theta
_{0})\right\Vert _{\mathbf{E}_{T}^{\nu }} &\leq &2^{-n}.
\end{eqnarray}%
Then, to conclude we just need to show that $\left( \phi
_{n}(\theta _{0})\right) _{n}$ is a Cauchy sequence in the Banach
space $\mathbf{z}_{\mathcal R}=L^{\infty }\left( \left[ 0,T\right]
,X_{\mathcal{R}}\right) $ endowed with its natural norm,
\begin{align*}
\left\Vert v\right\Vert _{\mathbf{Z}_{\mathcal
R}}=\sup_{0<t<\delta }\left( \left\Vert v(t)\right\Vert
_{X}+\left\Vert {\mathcal{R}^{\perp }}(v)(t)\right\Vert
_{X}\right).
\end{align*}
Firstly, using the Lemma \ref{lemma1} and the fact that $\left(
\phi _{n}(\theta _{0})\right) _{n}\in \mathbf{E}_{T}^{\nu },$ we
infer inductively that the sequence $\left( \phi _{n}(\theta
_{0})\right) _{n}$ belongs to the space $\mathbf{Z
}_{\mathcal R}.$
Secondly, once again the Lemma \ref{lemma1}, implies that the
sequence $\left( \omega _{n+1}\right) _{n}\equiv \left( \phi
_{n+1}(\theta _{0})-\phi_{n}(\theta _{0})\right) _{n}$ satisfies
the following inequality
\begin{eqnarray*}
\left\Vert \omega _{n+1}\right\Vert _{\mathbf{Z}_{\mathcal R}} &\leq &C\left(
\left\Vert \phi _{n}(\theta _{0})\right\Vert _{\mathbf{Z}_{\mathcal R}}+\left\Vert
\phi _{n-1}(\theta _{0})\right\Vert _{\mathbf{Z}_{\mathcal R}}\right) \left\Vert
\omega _{n}\right\Vert _{\mathbf{E}_{\delta }^{\nu }} \\
&&+C\left( \left\Vert \phi _{n}(\theta _{0})\right\Vert _{\mathbf{E}%
_{T}^{\nu }}+\left\Vert \phi _{n-1}(\theta _{0})\right\Vert _{\mathbf{E}%
_{\delta }^{\nu }}\right) \left\Vert \omega _{n}\right\Vert _{\mathbf{Z}%
_{\mathcal R}},
\end{eqnarray*}%
where $C=C(X,\alpha )>0.$
This inequality combined with the estimates \Ref{est20a}-\Ref{est20b} yields
\[
\left\Vert \omega _{n+1}\right\Vert _{\mathbf{Z}_{\mathcal R}}\leq C\left( \frac{1}{2}%
\right) ^{n}\left( \left\Vert \phi _{n}(\theta _{0})\right\Vert _{\mathbf{Z}%
_{\mathcal R}}+\left\Vert \phi _{n-1}(\theta _{0})\right\Vert _{\mathbf{Z}%
_{\mathcal R}}\right) +4C\mu \left\Vert \omega _{n}\right\Vert _{\mathbf{Z}%
_{\mathcal R}}
\]
Finally, if we choose $\mu > 0 $ such that $4C\mu <1$ one can
conclude the proof by using the following Lemma which is inspired
from \cite{furioli}.\fin
\begin{lemma}\label{pointfixe}
 Let $(x_{n})_{n}$ be a sequence in a normed vector space $\left(
Z,\left\Vert .\right\Vert \right) .$ If there exist a constant $\lambda \in
\lbrack 0,1[$ and $(\sigma _{n})_{n}\in l^{1}(\mathbb{N})$ such that:%
\begin{align}
\label{sequence}
~\forall n\in \mathbb{N}^{\ast },\quad
\left\Vert x_{n+1}-x_{n}\right\Vert \leq \sigma _{n}\left( \left\Vert
x_{n}\right\Vert +\left\Vert x_{n-1}\right\Vert \right) +\lambda \left\Vert
x_{n}-x_{n-1}\right\Vert,
\end{align}
then the series $\sum_{n}\left\Vert x_{n+1}-x_{n}\right\Vert $ converges. In
particular, $(x_{n})_{n}$ is a Cauchy sequence in $Z.$
\end{lemma}
{\bf Proof :} let us define the sequence $M_{n}=\sup_{k\leq n}\left\Vert
x_{k}\right\Vert .$ It follows inductively from \Ref{sequence}
\begin{eqnarray}
\nonumber
\left\Vert x_{n+1}-x_{n}\right\Vert &\leq &2\sum_{k=0}^{n-1}\sigma
_{n-k}M_{n-k}\lambda ^{k}, \\
\label{est21}
&\leq &\varpi _{n}M_{n},
\end{eqnarray}%
where $\varpi_{n}=2\sum_{k=0}^{n-1}\sigma _{n-k}\lambda ^{k}$.\\
Noticing that since $\left( \varpi _{n}\right) _{n}$ is a convolution of two
sequences in $l^{1}(\mathbb{N})$ then $\left( \varpi _{n}\right) _{n}$
belongs to $l^{1}(\mathbb{N}).$ Therefore,  we just need to show
that the sequence $\left( M_{n}\right) _{n}$ is bounded. This is somehow
obvious. In fact, using the triangular inequality $\left\Vert
x_{n+1}\right\Vert \leq \left\Vert x_{n}\right\Vert +\left\Vert
x_{n+1}-x_{n}\right\Vert ,$  \Ref{est21} yields%
\[
M_{n+1}\leq  ( 1+ \varpi_{n} ) M_{n}.
\]%
Which in turn implies%
\[
M_{n}\leq \Pi _{k=0}^{n-1}\left( 1+\varpi _{k}\right) \leq e^{\sum_{k\geq
0}\varpi _{n}}.
\]%
The proof is then achieved.\fin\\
Now let us see how the two previous lemmas allow to prove the Proposition \ref{propositionabstraite}.\\
{\bf Proof :} as usual we consider only the case of $\theta _{0}\in
X_{\mathcal{R}}$. Let $T > 0.$
By the Theorem \ref{theo1}, the solution $\theta $ is continous from $%
\mathbb{R}^{+}$ into $\mathbf{\tilde{B}}^{\alpha
},$ then $K \equiv \theta (\left[ 0,T%
\right] )$ is a compact subset of $\mathbf{\tilde{B}}^{\alpha }.$
Therefore, in view
of the Lemma \ref{compact}, there exists $\delta >0$ such that%
\begin{align}
\label{est22}
~\forall \tau \in \left[ 0,T\right] ,\quad
\left\Vert e^{-t\adelta}\theta (\tau )\right\Vert _{\E_{\delta
}^{\nu }}\leq \mu_0,
\end{align}
where $\mu_0 $ is the real given by Lemma \ref{localpersistence}. Now, we consider a repartition $%
0=t_{0}<\cdots <t_{N+1}=T$ of the interval $\left[ 0,T\right] $ such that $%
\sup_{i}t_{i+1}-t_{i}\leq \frac{\delta }{2}.$ We will show inductively that
\begin{align}
\label{est22a} \theta \in L^{\infty }\left( \left[
t_{i},t_{i+1}\right] ,X_{\mathcal{R}}\right) ,
\end{align}
which implies in turn the desired result $\theta \in L^{\infty
}\left( \left[ 0,T\right] ,X_{\mathcal{R}}\right).\;$ First, by
the Lemma \ref{localpersistence}, the claim  \Ref{est22a} is true for $i=0.$ Assume that, it is also true for $%
i\leq N.$ Then there exists $\tau _{0}$ in $\left]
t_{i},t_{i+1}\right[ $ such that $\tilde{\theta}_{0} \equiv \theta
(\tau _{0}) \in X\cap \mathbf{\tilde{B}}^{\alpha }.$ We notice
that $\tilde{\theta} \equiv  \theta(. + \tau_0) $ is the unique
solution given by Theorem \ref{theo1} of the Quasi-geostrophic
equation with initial data $\tilde{\theta}_{0}. $ Then according
to Lemma \ref{localpersistence} and \Ref{est22}, we get $ \theta
\in L^{\infty }\left( \left[ \tau_{0},\tau_{0} + \delta \right]
,X_{\mathcal{R}}\right)$. Hence, we are ready to conclude since $
\left[t_{i+1},t_{i+2}\right] \subset \left[\tau_{0}, \tau_{0} +
\delta \right].  $ \fin
\subsection{\bf Propagation of $B_p^{s,q}$ regularity for $ s <
0$}
\begin{proposition}\label{propsnegative}
Let $X $ be $ B_p^{s,q}$ or $ \dot{B}_p^{s,q}$  with $ s < 0 $ and
$ 1 \leq p, q \leq \infty.$ If $\theta_0 $ belongs to
$X_{\mathcal{R}}$ then the solution
$$
\theta \in \bigcap_{T>0} L^{\infty}([0,T], X_{\mathcal{R}}).
$$
\end{proposition}
As in the case $ s > 0, $ by using the compactness lemma \ref{compact} we just need to prove the following local persistency result :
\begin{lemma}
If $\theta_0  \in X_{\mathcal{R}}$ then there exists $\delta > 0 $
such that,
$$
\theta \in  L^{\infty}([0,\delta], X_{\mathcal{R}}).
$$
Moreover, the time $\delta $ is bounded below by,
\[
\sup \left\{ T>0/\left\Vert e^{-t\adelta }\theta _{0}\right\Vert
_{\E_{T}^{\nu }}\leq \mu_0 \right\},
\]%
where $\;\mu_0 \;$ is given by Lemma \ref{lemma3}.
\end{lemma}
{\bf Proof :} we consider only the case of $X =  B_p^{s,q}.$ The proof in the
other case is similar. Let  $ T > 0$ such that
$$
\left\Vert e^{-t\adelta }\theta _{0}\right\Vert _{\E_{T}^{\nu
}}\leq \mu_0.
$$
According to the Lemma \ref{lemma3} the sequence $\;(\phi_n(\theta_0))_n\;$ satisfies
\begin{align}
\label{est35}
||\phi_{n+1}(\theta_0)- \phi_n(\theta_0) ||_{\E_{T}^{\nu}} \> \leq \frac{1}{2^n},
\end{align}
and converges to the solution $\theta$ in $\;\E_{T}^{\nu}.\;$ Our first task is to prove that $(\phi_n(\theta_0))_n$ is a Cauchy sequence in the space,
$$
X_{\sigma,p}^T = \{  v  : ( 0 , T ] \rightarrow L^p ; \;
||v||_{X_{\sigma,p}^T } \equiv \sup_{0<t<T} t^{\frac{\sigma}{2\alpha}}
( ||v(t)||_p + ||\mathcal{R}^{\perp } (v)(t) ||_p ) < \infty  \},
$$
where $\sigma = -s.\;$\\
Thanks to the Besov  characterization \Ref{car02} and  Lemma \ref{lemma1}, we can show inductively that $(\phi_n(\theta_0))$ belongs to  $
X_{\sigma,p}^T$ and satisfies,
\begin{align}
\label{est35a}
||\phi_{n+1}(\theta_0) - \phi_n(\theta_0) ||_{X_{\sigma,p}^T}  \leq C ||\phi_n(\theta_0)  - \phi_{n-1}(\theta_0) ||_{\E_{T}^{\nu}}
\max( ||\phi_n(\theta_0)  ||_{X_{\sigma,p}^T} , || \phi_{n-1}(\theta_0) ||_{X_{\sigma,p}^T}).
\end{align}
Thus, By  \Ref{est35} and   Lemma \ref{pointfixe} we deduce that
$(\phi_n(\theta_0) )_n$ is a Cauchy sequence in $X_{\sigma,p}^T.$ Therefore
its limit $\theta \in X_{\sigma,p}^T.$ Now by a simple computation
using the characterization \Ref{car02} we deduce that $ \theta \in
L^{\infty}( [ 0, T_0], (B_p^{s,\infty})_{\mathcal{R}}). $
Moreover, for $\epsilon
> 0$ such that 
\bea
\label{eps}
 -1 < s \pm \epsilon < 0,
 \eea 
one can show that the nonlinear part $N(\theta)(t) = {\mathcal{B}}_{\alpha}%
\left[ \theta ,\theta \right] (t)$ satisfies
\begin{align}
\label{matin}
|| N(\theta)(t) ||_{B_p^{s\pm\epsilon},\infty}  + || {\mathcal R}^{\perp
} N(\theta)(t) ||_{B_p^{s\pm\epsilon},\infty} \leq  C_{s,\epsilon}
t^{-\pm\frac{\epsilon}{2\alpha}} ||\theta||_{\E_{T}^{\nu}}
||\theta||_{X_{\sigma,p}^T}.
\end{align}
Indeed, we have $\tau \in ] 0, 1 [$
\bea
\nonumber
\tau^{-\frac{s\pm\epsilon}{2\alpha}} || e^{-\tau\adelta} N(\theta)(t)||_p & \leq & C \int_0^t (t + \tau - r )^{-\frac{1}{2\alpha}}  \tau^{-\frac{s\pm\epsilon}{2\alpha}}  r^{-\nu}  r^{-\frac{\sigma}{2\alpha}} dr \; ||\theta||_{\E_{T}^{\nu}}
||\theta||_{X_{\sigma,p}^T}, \\
\nonumber
 &\leq & C \int_0^t (\frac{\tau}{t+\tau -r})^{-\frac{s\pm\epsilon}{2\alpha}} 
 ( t + \tau -r)^{\frac{-1 -(s\pm\epsilon)}{2\alpha}} r^{-\nu}  r^{-\frac{\sigma}{2\alpha}} dr \; ||\theta||_{\E_{T}^{\nu}}
||\theta||_{X_{\sigma,p}^T}, \\
\label{m3}
&\leq & C \int_0 ^t   ( t  -r)^{\frac{-1 - ( s \pm\epsilon)}{2\alpha}} r^{-\nu}  r^{-\frac{\sigma}{2\alpha}} dr \; ||\theta||_{\E_{T}^{\nu}}
||\theta||_{X_{\sigma,p}^T}, \\
\label{m4}
&\leq & C  t ^{-\frac{\pm\epsilon}{2\alpha}}  \; ||\theta||_{\E_{T}^{\nu}}||\theta||_{X_{\sigma,p}^T},
\eea
Where to obtain \Ref{m3}, we have used the facts that , $ 0 \leq \frac{\tau}{t+\tau -r} \leq 1,\quad  t + \tau - r \geq t - r$ and  \Ref{eps}. Similarly, we have the same estimate \Ref{m4} for the ${\mathcal R}^{\perp
} N(\theta)(t).$ Hence, by Proposition \ref{p} we get \Ref{matin}.
Thus, by interpolation we obtain $N(\theta) \in L^{\infty}([ 0, T],
(B_p^{s,1})_{\mathcal R})$ which implies $\theta \in L ^{\infty}([
0, T], (B_p^{s,q})_{\mathcal R}).$ 
\fin 

\subsection{\bf The case of null regularity $ s = 0$}
In this subsection we aim to prove the following result,
\begin{proposition}
Let $X $ be $ B_p^{0,q}$ or $ \dot{B}_p^{0,q}$  with  $ 1 \leq p,
q \leq \infty.$ If $\theta_0\in  X$ then the solution
$$
\theta \in \bigcap_{T>0} L^{\infty}([0,T], X).
$$
\end{proposition}

Thanks to the following imbedding,
$$
\dot{B}_p^{0,1} \subset \dot{B}_p^{0,q}  \subset \dot{B}_p^{0,\infty},
$$
and
$$
 \dot{B}_p^{0,1} \subset B_p^{0,q}  \subset \dot{B}_p^{0,\infty},
$$
the proof of the above proposition is an immediate consequence of the following lemma,
\begin{lemma}
If  $\; \theta_0 \in \dot{B}_p^{0,\infty}\;$ then $\; N(\theta)= {\mathcal{B}}_{\alpha}%
\left[ \theta ,\theta \right] (t) \;$ belongs to $
\bigcap_{T>0}L^{\infty}([0,T], \dot{B}_p^{0,1} ).$
\end{lemma}
{\bf Proof :}
By using the Young inequality we deduce that
\begin{align*}
\dot{B}_{p}^{0,\infty} \cap {\dot{B}}^{-(2\alpha-1),\infty }_{\infty} \subset  \dot{B}_{2p}^{\frac 12 -\alpha,\infty}.
\end{align*}
Observe that $s^{*} =  \frac 12 -\alpha < 0 $ and hence according to the proof of the proposition \ref{propsnegative}
 and to the continuity of the Riesz transforms on homogeneous Besov spaces, we have $\theta \in \bigcap_{T > 0} X_{\sigma^{*},2p}^T\;$ where $\; \sigma^{*} = \alpha - \frac 12.$
Let $ \;T > 0 \;$ and    $ 0 < \sigma < 2\alpha-1.$ The basic estimate,
\begin{align*}
|| \sqrt{-\Delta}^{\pm\sigma}\; \nabla e^{-t{\adelta}} f ||_p &  \leq  C_0 t^{-\frac{\pm\sigma + 1}{2\alpha}} || f ||_p.
\end{align*}
yields immediately
\begin{align*}
|| (\sqrt{-\Delta})^{\pm \sigma} N(\theta)(t) ||_p  &
\leq  C
  t^{-\frac{\pm\sigma}{2\alpha}} || \theta ||_{X_{\sigma^*,2p}^T}^2.
\end{align*}
Now, we use the interpolation result (see Theorem 6.3 in
\cite{bergh})
\begin{align*}
[
(\sqrt{-\Delta})^{\sigma} L^p,
(\sqrt{-\Delta})^{-\sigma} L^p
]_{\frac{1}{2},1} =  \dot{B}_{p}^{0,1},
\end{align*}
to deduce,
\begin{align}
\label{nonl2} \forall \; 0 < t < T, \quad || N(\theta)(t) ||_{
\dot{B}_{p}^{0,1}} & \leq  C  || \theta ||_{X_{\sigma^*,2p}^T}^2,
\end{align}
that implies,
\begin{align}
\label{nonlN}
N(\theta) \in L^{\infty}([0,T], \dot{B}_{p}^{0,1} ).
\end{align}
\fin 
\begin{remark} As in the context of the Navier-Stokes equations \cite{cannone2000}, we observe thanks to \Ref{nonlN} and \Ref{nonl2} that in the case $\; -1 < s \leq 0\; $, the fluctuation term $w(t)$ is more regular than the tendency  $e^{-t\adelta}\theta_0.\;$
\end{remark}
\section{Proof of Theorem \ref{theo3}}\label{sectlp}
The existence part is a direct consequence of Theorem \ref{theo1},
Theorem \ref{theo2} and the following embedding (consequence of
Bernstein's inequality and the boundedness of the Riesz transforms
on Lebesgue's and Sobolev's spaces)
\begin{eqnarray*}
L^{p}(\mathbb{R}^{2}) &\subset &\mathbf{\tilde{B}}_{\alpha
}~\forall \; p\geq
p_{c}, \\
H^{s}(\mathbb{R}^{2}) &=&B_{2}^{s,2}(\mathbb{R}^{2})\subset \mathbf{\tilde{B}%
}_{\alpha }~\forall \; s\geq s_{c}.
\end{eqnarray*}%
Let us establish the uniqueness part. First we notice that since
for $s\geq s_{c}$ then
\[
H^{s}(\mathbb{R}^{2})\hookrightarrow H^{s_{c}}(\mathbb{R}^{2})%
\hookrightarrow L^{p_{c}}(\mathbb{R}^{2}).
\]%
therefore, we just need to prove the uniqueness in the spaces
$\left( C([0,T],L^{p}(\mathbb{R}^{2}))\right) _{p\geq p_{c}}.$
This will be deduced
from the following continuity result of the bilinear operator $\mathcal{B}%
_{\alpha }.$

\begin{lemma}
Let $p\in ]p_{c},\infty \lbrack ,q\in ]1,\infty \lbrack $ and
$T>0.$ There exists a constant $C$ independent of $T$ such that:
\end{lemma}

\begin{itemize}
\item for any $u,v$ in $L_{T}^{\infty }L^{p},$%
\begin{equation}
\left\Vert \mathcal{B}_{\alpha }[u,v]\right\Vert _{L_{T}^{\infty
}L^{p}}\leq C\; T^{\sigma }\left\Vert u\right\Vert _{L_{T}^{\infty
}L^{p}}\left\Vert v\right\Vert _{L_{T}^{\infty }L^{p}}\text{,}
\label{ESS}
\end{equation}%
where $\sigma =\frac{1}{\alpha }(\frac{1}{p_{c}}-\frac{1}{p}).$

\item for any $u,v$ in $L_{T}^{\infty }L^{p_{c}},$%
\begin{equation}
\left\Vert \mathcal{B}_{\alpha }[u,v]\right\Vert
_{L_{T}^{q}L^{p_{c}}}+\left\Vert \mathcal{B}_{\alpha
}[v,u]\right\Vert _{L_{T}^{q}L^{p_{c}}}\leq C\left\Vert
u\right\Vert _{L_{T}^{\infty }L^{p_{c}}}\left\Vert v\right\Vert
_{L_{T}^{q}L^{p_{c}}.}  \label{ESS2}
\end{equation}

\item for any $u\in L_{T}^{\infty }L_{\mathcal{R}}^{\infty }$ and
$v\in
L_{T}^{q}L^{p_{c}},$%
\begin{equation}
\left\Vert \mathcal{B}_{\alpha }[u,v]\right\Vert
_{L_{T}^{q}L^{p_{c}}}+\left\Vert \mathcal{B}_{\alpha
}[v,u]\right\Vert _{L_{T}^{q}L^{p_{c}}}\leq C\;
T^{1-\frac{1}{2\alpha }}\left\Vert u\right\Vert _{L_{T}^{\infty
}L_{\mathcal{R}}^{\infty }}\left\Vert v\right\Vert
_{L_{T}^{q}L^{p_{c}}.}  \label{ESS3}
\end{equation}
\end{itemize}

{\bf Proof :}
Estimate (\ref{ESS}) follows easily from the continuity of the
Riesz transforms on the Lebesgue spaces $L^{r}(\mathbb{R}^{2})$
with $1<r<\infty ,$ the Young and the H{\"o}lder inequality and the
estimate \ref{est1a} on the $L^{r}(\mathbb{R}^{2})$ norm of the
kernel of the operator $\nabla e^{-(t-s)\left( -\Delta \right)
^{\alpha }} .$ Estimate (\ref{ESS2}) is a consequence of the
continuity of the Riesz transforms on the space
$L^{p_{c}}(\mathbb{R}^{2})$, the H{\"o}lder
inequality, the Sobolev embedding%
\[
\left\Vert \frac{\nabla }{\left( -\Delta \right) ^{\alpha
}}f\right\Vert _{p_{c}}\lesssim \left\Vert f\right\Vert
_{\frac{p_{c}}{2}}
\]%
and the following maximal regularity property of the operator
$\left( -\Delta \right) ^{\alpha }$
\[
\left\Vert \int_{0}^{t}\left( -\Delta \right) ^{\alpha
}e^{-(t-s)\left( -\Delta \right) ^{\alpha }}vds\right\Vert
_{L_{T}^{q}L^{p_{c}}}\lesssim \left\Vert v\right\Vert
_{L_{T}^{q}L^{p_{c}}}
\]%
which can be proved by following the proof of Theorem 7.3 in
\cite{lemarie}. Let us now
prove estimate (\ref{ESS3}). For any $t\in \lbrack 0,T]$ we have%
\begin{eqnarray*}
\left\Vert \mathcal{B}_{\alpha }[u,v](t)\right\Vert _{L^{p_{c}}}
&\lesssim &\int_{0}^{t}\frac{1}{(t-s)^{1/2\alpha }}\left\Vert
\mathcal{R}^{\perp
}(u)(s)\right\Vert _{\infty } \left\Vert v(s)\right\Vert _{p_{c}}ds \\
&\lesssim &\left\Vert \mathcal{R}^{\perp }(u)\right\Vert
_{L_{T}^{\infty }L^{\infty }}\left( 1_{[0,T]}s^{-\frac{1}{2\alpha
}}\right) \ast \left( 1_{[0,T]}  \left\Vert v(s)\right\Vert _{p_{c}}\right) (t)
\end{eqnarray*}%
where the star $\ast $ denotes the convolution in $\mathbb{R}.$
Hence the
Young inequality yields%
\[
\left\Vert \mathcal{B}_{\alpha }[u,v]\right\Vert
_{L_{T}^{q}L^{p_{c}}}\lesssim \left\Vert \mathcal{R}^{\perp
}(u)\right\Vert _{L_{T}^{\infty }L^{\infty }}T^{1-\frac{1}{2\alpha
}}\left\Vert v\right\Vert _{L_{T}^{q}L^{p_{c}}}.
\]%
Similarly, we obtain%
\begin{eqnarray*}
\left\Vert \mathcal{B}_{\alpha }[v,u]\right\Vert
_{L_{T}^{q}L^{p_{c}}} &\lesssim &T^{1-\frac{1}{2\alpha
}}\left\Vert u\right\Vert _{L_{T}^{\infty }L^{\infty }}\left\Vert
\mathcal{R}^{\perp }\left( v\right) \right\Vert
_{L_{T}^{q}L^{p_{c}}} \\
&\lesssim &T^{1-\frac{1}{2\alpha }}\left\Vert u\right\Vert
_{L_{T}^{\infty }L^{\infty }} \left\Vert
v\right\Vert _{L_{T}^{q}L^{p_{c}}}.
\end{eqnarray*}%
Estimate (\ref{ESS3}) is then proved.
\fin\\

Now we are ready to finish the proof of the uniqueness. Let $p\geq
p_{c}$ and $T>0$ be two reals number and let $\theta _{1}$ and
$\theta _{2}$ be two mild solutions
of the equation $(QG_{\alpha })$ with the same data $\theta _{0}$ such that $%
\theta _{1},\theta _{2}\in C([0,T],L^{p}(\mathbb{R}^{2})).$ We aim
to show that $\theta _{1}=\theta _{2}$ on $[0,T].$ For this, we
will argue by
contradiction. Then we suppose that $t_{\ast }<T$ where%
\[
t_{\ast }\equiv \sup \{t\in \lbrack 0,T];\forall s\in \lbrack
0,t],~\theta _{1}(s)=\theta _{2}(s)\}.
\]%
To conclude, we need to prove that there exists $\delta \in
]0,T-t_{\ast }]$
such that $\tilde{\theta}_{1}=\tilde{\theta}_{2}$ on $[0,\delta ],$ where $%
\tilde{\theta}_{1}$ and $\tilde{\theta}_{2}$ are the functions defined on $%
[0,T-t_{\ast }]$ by
\[
\tilde{\theta}_{1}(t)=\theta _{1}(t+t_{\ast
}),~\tilde{\theta}_{2}(t)=\theta _{2}(t+t_{\ast }).
\]%
We deal separately with the sub-critical case and the critical
case:

\textbf{The first case:} $p>p_{c}.$ Thanks to the continuity of
$\theta _{1}$ and $\theta _{2}$ on $[0,T],$ we have $\theta
_{1}(\tau _{\ast })=\theta
_{2}(t_{\ast }).$ Hence, the functions $\tilde{\theta}_{1}$ and $\tilde{%
\theta}_{2}$ are two mild solutions on $[0,\delta _{0}\equiv
T-t_{\ast }]$ of the equation $(QG_{\alpha })$ with the same data
$\theta _{1}(\tau _{\ast
}).$ Therefore, the function $\tilde{\theta}\equiv \tilde{\theta}_{1}-\tilde{%
\theta}_{2}$ satisfies the equation%
\begin{equation}
\tilde{\theta}=\mathcal{B}_{\alpha }[\tilde{\theta}_{1},\tilde{\theta}]-%
\mathcal{B}_{\alpha }[\tilde{\theta},\tilde{\theta}_{2}].
\label{equa}
\end{equation}%
Thus, according to (\ref{ESS}) we have for any $\delta \in ]0,\delta _{0}]$%
\begin{eqnarray*}
\left\Vert \tilde{\theta}\right\Vert _{L_{\delta }^{\infty }L^{p}}
&\leq &C\delta ^{\sigma }\left( \left\Vert
\tilde{\theta}_{1}\right\Vert _{L_{\delta }^{\infty
}L^{p}}+\left\Vert \tilde{\theta}_{2}\right\Vert _{L_{\delta
}^{\infty }L^{p}}\right) \left\Vert \tilde{\theta}\right\Vert
_{L_{\delta }^{\infty }L^{p}} \\
&\leq &C\delta ^{\sigma }\left( \left\Vert \theta _{1}\right\Vert
_{L_{T}^{\infty }L^{p}}+\left\Vert \theta _{2}\right\Vert
_{L_{T}^{\infty }L^{p}}\right) \left\Vert
\tilde{\theta}\right\Vert _{L_{\delta }^{\infty }L^{p}}
\end{eqnarray*}%
where $C>0$ is independent on $\delta .$

Consequently, for $\delta $ small enough, $\tilde{\theta}=0$ on
$[0,\delta ]$ which ends the proof in the sub-critical case.

\textbf{The second case:} $p=p_{c}.$ Choose a fix real $q>1$ and let $%
\varepsilon >0$ to be chosen later. By density of smooth functions in the space $C([0,T],L^{p_c}(\mathbb{R}^{2}))$, one can decompose $\tilde{%
\theta}_{1}$ and $\tilde{%
\theta}_{2}$ into $\tilde{%
\theta}_{1}=u_{1}+v_{1}$ and $\tilde{\theta}_{2}=u_{2}+v_{2}$ with%
\begin{eqnarray}
\left\Vert u_{1}\right\Vert _{L_{\delta _{0}}^{\infty
}L^{p_{c}}}+\left\Vert u_{2}\right\Vert _{L_{\delta _{0}}^{\infty
}L^{p_{c}}} &\leq &\varepsilon ,
\label{es1} \\
\left\Vert v_{1}\right\Vert _{L_{\delta _{0}}^{\infty }L_{\mathcal{R}%
}^{\infty }}+\left\Vert v_{2}\right\Vert _{L_{\delta _{0}}^{\infty }L_{%
\mathcal{R}}^{\infty }} &\equiv &\mathcal{M}<\infty .  \label{es2}
\end{eqnarray}%
As in the previous case, the function $\tilde{\theta}\equiv \tilde{\theta}%
_{1}-\tilde{\theta}_{2}$ satisfies the equations%
\begin{eqnarray*}
\tilde{\theta} &=&\mathcal{B}_{\alpha }[\tilde{\theta}_{1},\tilde{\theta}]+%
\mathcal{B}_{\alpha }[\tilde{\theta},\tilde{\theta}_{2}] \\
&=&\mathcal{B}_{\alpha }[u_{1},\tilde{\theta}]+\mathcal{B}_{\alpha }[\tilde{%
\theta},u_{2}]+\mathcal{B}_{\alpha }[v_{1},\tilde{\theta}]+\mathcal{B}%
_{\alpha }[\tilde{\theta},v_{2}].
\end{eqnarray*}%
Now by applying (\ref{ESS2})-(\ref{ESS3}) and using
(\ref{es1})-(\ref{es2})
we get, for any $\delta \in ]0,\delta _{0}],$ the following estimate%
\[
\left\Vert \tilde{\theta}\right\Vert _{L_{\delta }^{q}L^{p}}\leq
C\left( \varepsilon +\delta ^{1-\frac{1}{2\alpha
}}\mathcal{M}\right) \left\Vert \tilde{\theta}\right\Vert
_{L_{\delta }^{q}L^{p}}
\]%
where $C>0$ is a constant depending only on $\alpha ,p$ and $q.$

Thus, by choosing $\varepsilon $ small enough, we conclude that
there exists $\delta \in ]0,\delta _{0}]$ such that $\left\Vert
\tilde{\theta}\right\Vert
_{L_{\delta }^{q}L^{p}}=0$, which implies that $\tilde{\theta}_{1}=\tilde{%
\theta}_{2}$ on $[0,\delta ]$. The proof is then achieved.
\begin{remark} The idea of the proof of the uniqueness in the
critical case is inspired from the paper \cite{monniaux} of S.
Monniaux.
\end{remark}


\end{document}